\documentclass[12pt]{article}
\usepackage{hyperref}
\usepackage{cmap}                        
\usepackage[cp1251]{inputenc}            
\usepackage[english]{babel}
\usepackage[left=1in,right=1in,top=1in,bottom=1.5in]{geometry} 
\usepackage{amssymb,amsmath, amsthm, amscd,ifthen}
\usepackage{graphicx}

\newcommand{\R}{{\mathbb R}}

\newtheorem{thm}{Theorem}
\newtheorem{lem}{Lemma}
\newtheorem{cor}{Corollary}
\date{}

\title{From Tarski's plank problem to simultaneous approximation}

\author{Andrey B. Kupavskii\thanks{EPFL, Lausanne and MIPT, Moscow.   Supported in part by the grant N 15-01-03530 of the Russian Foundation for Basic Research. E-mail: {\tt kupavskii@ya.ru}.} \and J\'anos Pach\thanks{EPFL, Lausanne and R\'enyi Institute, Budapest. Supported
by Hungarian Science Foundation EuroGIGA Grant OTKA NN 102029, by Swiss National Science Foundation Grants 200020-144531 and 200021-137574. E-mail:
{\tt pach@cims.nyu.edu}.} }
\date{}

\begin{document}
\maketitle

\begin{abstract} \noindent
    A {\em slab} (or plank) of width $w$ is a part of the $d$-dimensional space that lies between two parallel hyperplanes at distance $w$ from each other. It is conjectured that any slabs $S_1, S_2,\ldots$ whose total width is divergent have suitable translates that altogether cover $\mathbb{R}^d$. We show that this statement is true if the widths of the slabs, $w_1, w_2,\ldots$, satisfy the slightly stronger  condition $\limsup_{n\rightarrow\infty}\frac{w_1+w_2+\ldots+w_n}{\log(1/w_n)}>0$. This can be regarded as a converse of Bang's theorem, better known as Tarski's plank problem.

    We apply our results to a problem on simultaneous approximation of polynomials. Given a positive integer $d$, we say that a sequence of positive numbers $x_1\le x_2\le\ldots$ {\em controls} all polynomials of degree at most $d$ if there exist $y_1, y_2,\ldots\in\mathbb{R}$ such that for every polynomial $p$ of degree at most $d$, there exists an index $i$ with $|p(x_i)-y_i|\leq 1.$ We prove that a sequence has this property if and only if $\sum_{i=1}^{\infty}\frac{1}{x_i^d}$ is divergent. This settles an old conjecture of Makai and Pach.
\end{abstract}

\section{Tarski's plank problem and its affine version}
The closed set of points $S$ lying between two parallel hyperplanes in $\mathbb{R}^d$ at distance $w$ from each other is called a {\em slab} (or {\em plank}) of {\em width} $w$. Given a convex body $C\subset\mathbb{R}^d$, its {\em width} $w(C)$ is the smallest number $w$ such that there is a slab of width $w$ that covers $C$.
\medskip

In 1932, {\em Alfred Tarski}~\cite{Ta} made the following attractive conjecture.
\medskip

\noindent{\bf Tarski's plank problem.} {\em If a sequence of slabs covers a convex body $C$, then the total width of the slabs is at least $w(C)$.}

\medskip

Tarski found a beautiful proof of his conjecture in the special case where $C$ is a closed disk. His argument was based on the following fact discovered by {\em Archimedes} more than two thousand years ago. Let $D$ be a disk of radius $\frac12$ in the $(x,y)$-plane $\mathbb{R}^2$, and let $H$ denote the hemisphere of radius $\frac12$, concentric with $D$, that lies in the closed half-space above $\mathbb{R}^2$. No matter how we place a $2$-dimensional slab $S$ of width $w$ in the plane such that both of its boundary lines intersect $D$, the surface area of the vertical projection of $S\cap D$ to the hemisphere $H$ is always equal to $\frac12\pi w$; see Fig. 1. Suppose now that a system of slabs $S_1, S_2,\ldots$ completely covers $D$. We can assume without loss of generality that the boundary lines of every slab intersect $D$, otherwise, we could replace one of the slabs with a narrower one. The vertical projections of the slabs cover $H$, therefore the sum of the areas of these projections, $\frac12{\pi w(S_1)}+\frac12{\pi w(S_2)}+\ldots$, is at least $Area(H)=\frac{\pi}2$. Thus, we have that $$w(S_1)+w(S_2)+\ldots\ge 1=w(D),$$ as required.
\medskip

\begin{figure}[h]
\centering
\smallskip
\includegraphics[height=3in]{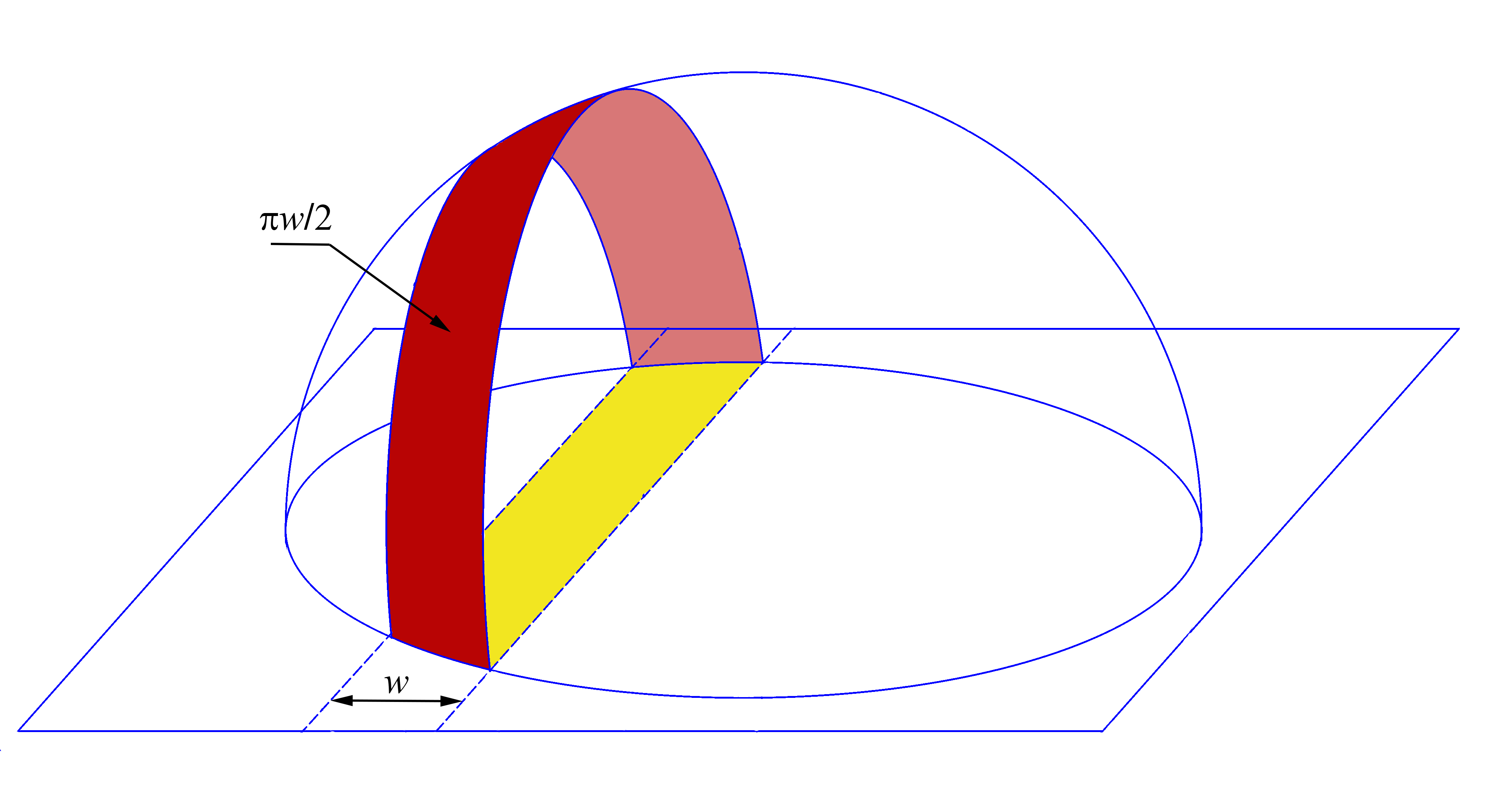}
\smallskip
\caption{Tarski's proof for a disk of unit diameter.}
\medskip
\end{figure}

A slight modification of the above argument also settles the analogous problem in the case where $C$ is a $3$-dimensional ball $B^3$ of unit diameter. Now let $H$ denote the whole sphere bounding $B^3$, so that $Area(H)=\pi$. Suppose that $S_1, S_2,\ldots$ is a system of $3$-dimensional slabs that altogether cover $B^3$, and they have the property that the boundary planes of each $S_i$ intersect $B^3$. By the Archemidean observation, for every $i$, the surface area of $S_i\cap H$ is equal to $\pi w(S_i)$. Since the slabs must also cover $H$, we have $$\sum_i Area(S_i\cap H)=\sum_i \pi w(S_i)\ge Area(H)=\pi,$$ which again yields that $\sum_i w(S_i)\ge w(B^3)=1.$

\medskip

One may naively hope that the above argument allows a straightforward generalization to higher dimensional balls or even to convex bodies of other shapes. This is not the case. It took almost twenty years before {\em Th{\o}ger Bang}~\cite{Ba1}, \cite{Ba2} managed to prove Tarski's conjecture, using quite different ideas. He also formulated a more general conjecture that he was unable to settle. To state Bang's conjecture, we need some notation.
\medskip

A unit vector $\mathbf{v}$ perpendicular to the bounding hyperplanes of a slab $S$ is called the {\em normal vector} of $S$. Given a convex body $C$ and a unit vector $\mathbf{v}\in \mathbb{R}^d$, let $w(C,\mathbf{v})$ denote the {\em width} of $C$ {\em in direction} $\mathbf{v}$, that is, the smallest number $w$ such that there is a slab of width $w$ with normal vector $\mathbf{v}$ that contains $C$. For example, if $C$ is a $d$-dimensional ball of unit diameter, then $w(C,\mathbf{v})=1$ for any unit vector $\mathbf{v}\in \mathbb{R}^d$.
\medskip

\noindent{\bf Bang's affine plank problem.} {\em Let $S_1, S_2, \ldots$ be a sequence of slabs in $\mathbb{R}^d$ with normal vectors ${\mathbf{v}}_1, {\mathbf{v}}_2, \ldots$ and widths $w(S_1), w(S_2),\ldots,$ respectively, and suppose that their union covers a convex body $C$.

Then we have
$$\sum_i \frac{w(S_i)}{w(C,{\mathbf{v}_i})}\ge 1.$$}
\medskip

This conjecture is still open. Bang~\cite{Ba3} verified it for systems consisting of only {\em two} slabs (see also ~\cite{Mo}), but it is not known to be true even for triples. It was a sensational breakthrough, when in 1991 {\em Keith Ball}~\cite{B} settled the conjecture in the affirmative for centrally symmetric convex bodies $C$. On the other hand, some results of {\em Richard Gardner}~\cite{Ga} indicate that the affine plank problem cannot be solved by any argument similar to Tarski's.

\section{A converse of Tarski's problem}

Tarski's conjecture, that is, Bang's theorem, states that if a sequence of slabs cover a convex body $C\subset\mathbb{R}^d$, then their total width must be large. One can reverse this question, as follows. Suppose that we have a sequence of slabs $S_1, S_2,\ldots$ of large total width. Is it always possible to cover $C$ with congruent copies of the $S_i$? The answer is simple. Take a slab $S$ of width $w(C)$ that contains $C$, and denote its normal vector by $\mathbf{v}$. If $\sum_i w(S_i)\ge w(C)$, we can rotate each $S_i$ into a position perpendicular to $\mathbf{v}$. Now we can translate these slabs so that their union will cover $S$ and, hence $C$.
\medskip

The problem becomes more interesting if we permit only translations, but no rotation. We say that a sequence of slabs $S_1, S_2,\ldots$ permits a {\em translative covering} of a subset $C\subseteq\R^d$ if there are suitable translates $S_i'$ of $S_i\; (i=1,2,\ldots)$ such that $\cup_{i}S_i'\supseteq C$.
\medskip

If each of the slabs has width greater than some positive constant $\epsilon$, then each can be used to cover a full-dimensional ball of diameter $\epsilon$. Since any convex body $C$ can be covered by finitely many balls of diameter $\epsilon$, these slabs permit a translative covering of $C$.

However, it is not clear whether the condition that the total width of the slabs is {\em divergent} is sufficient to guarantee that they permit a translative covering of a ball of unit diameter. If this is true, then any sequence of slabs with divergent total width also permits a translative covering of the {\em whole space} $\mathbb{R}^d$. Indeed, any such sequence can be partitioned into infinitely many subsequences, each having divergent total width. Choose any covering of $\mathbb{R}^d$ with balls $B_1, B_2, \ldots$ of unit diameter, and for each $B_i$, use translates of the slabs belonging to the $i$th subsequence to cover $B_i\; (i=1,2,\ldots)$. {\em Endre Makai} and {\em J\'anos Pach}~\cite{MP} made the following conjecture; see also ~\cite{BMP}, Section 3.4.
\medskip

\noindent{\bf Makai--Pach translative plank conjecture.}
{\em Let $d$ be a positive integer. A sequence of slabs in $\mathbb{R}^d$ with widths $w_1, w_2,\ldots$  permits a translative covering of $\mathbb{R}^d$ if and only if $\sum_{i=1}^{\infty}w_i=\infty$.}
\medskip

The ``only if'' part is quite easy. It also follows directly from Bang's theorem. If $\sum_{i=1}^{\infty}w_i<D$, then the slabs do not even permit a translative covering of a ball of diameter $D$.

As for the ``if'' part, in the $2$-dimensional case it was proved in~\cite{MP} and, according to~\cite{Gr0}, independently, by {\em Paul Erd\H os} and {\em Ernst G. Straus} (unpublished). In this case, there is a constant $c>0$ such that any system of slabs in the plane with total width at least $c$ permits a translative covering of a disk of unit diameter, and the conjecture is true. (See~\cite{Gr1,Gr2} for some refinements.) For $d\ge 3$, the best known result is due to {\em Helmut Groemer}~\cite{Gr0}. He proved that for any $d\ge 3,$ any sequence of slabs of widths $w_1, w_2, \ldots$, satisfying  $$\sum_{i=1}^{\infty}w_i^{\frac{d+1}2}=\infty,$$ permits a translative covering of $\mathbb{R}^d$. In particular, for $d=3$, any sequence of slabs with widths $w_i=\frac1{\sqrt{i}}\; (i=1,2,\ldots)$ permits a translative covering of $\mathbb{R}^3$. Our next result shows that the same is true for much narrower slabs, for example, for slabs of widths $w_i=\frac1{i}$ for every $i$. This comes rather close to the truth: the last statement is false, e.g., for the sequence $w_i=1/i^{1+\varepsilon}\; (i=1,2,\ldots)$ with any $\varepsilon>0,$ because then we have $\sum_{i=1}^{\infty}w_i<\infty$.

\begin{thm}\label{th1}
Let $d$ be a positive integer, and let $w_1\ge w_2\ge\ldots$ be a monotone decreasing infinite sequence of positive numbers such that $$\limsup_{n\rightarrow\infty}\frac{w_1+w_2+\ldots+w_n}{\log(1/w_n)}> 0.$$
Then any sequence of slabs $S_i$ of width $w_i\; (i=1,2,\ldots)$ permits a translative covering of $\mathbb{R}^d$.
\end{thm}
\noindent Here and in what follows $\log$ stands for the natural logarithm.
\smallskip

Our proof is based on the following statement.

\begin{thm}\label{th2} Let $d$ be a positive integer, and let $w_1\ge w_2\ge\ldots\ge w_n$ be positive numbers such that $$w_1+w_2+\ldots+w_n\ge 3d\log(2/w_n).$$
Then any sequence of slabs $S_1,\ldots, S_n\subset\mathbb{R}^d$ with widths $w_1,\ldots,w_n$, resp., permits a translative covering of a $d$-dimensional ball of diameter $1-w_n/2$.
\end{thm}

Theorems 1 and 2 will be established in Section 4.

\section{Application to approximation of polynomials}

For a fixed positive integer $d$, let ${\mathcal P}_d$ denote the class of polynomials of degree at most $d$. Following~\cite{MP}, we say that a sequence of positive numbers $x_1, x_2,\ldots$ is {\em ${\mathcal P}_d$-controlling} if there exist reals $y_1, y_2,\ldots$ with the property that for every polynomial $p\in{\mathcal P}_d$, one can find an $i$ with $$|p(x_i)-y_i|\le 1.$$
Roughly speaking, this means that the graph of every polynomial $p\in{\mathcal P}_d$ comes vertically close to at least one point of the set $\{(x_i,y_i)\,:\,i=1,2\ldots\}$.
\medskip

In Section 5, we study the following question. How {\em sparse} can a ${\mathcal P}_d$-controlling sequence be? A similar question, motivated by a problem of {\em L\'aszl\'o Fejes T\'oth}~\cite{FT}, was studied in~\cite{EP}. We will see that this question is intimately related to the translative plank conjecture discussed in the previous section.

First, we show how the following assertion can be deduced from Theorem~\ref{th1}.

\begin{cor}\label{th3} Let $d$ be a positive integer and $x_1\le x_2\le \ldots$ be a monotone increasing infinite sequence of positive numbers. If $$\limsup_{n\rightarrow\infty}\left(\frac1{x_1^{d}}+\frac1{x_2^{d}}+\ldots+\frac1{x_n^{d}}\right)\Big/{\log x_n}>0,$$
then the sequence $x_1, x_2,\ldots$ is {\em ${\mathcal P}_d$-controlling}.
\end{cor}

\medskip
\noindent{\bf Proof of Corollary~\ref{th3}.}
Let $x_1\le x_2\le\ldots$ be an infinite sequence of positive numbers satisfying the assumptions. We have to find a sequence of reals $y_1, y_2,\ldots$ such that for any polynomial $p(x) = \sum_{j=0}^da_jx^j$ with real coefficients $a_j$, there exists a positive integer $i$ with $|p(x_i)-y_i|\le 1$. See Fig. 2.

\begin{figure}[h]
\centering
\smallskip
\includegraphics[height=3in]{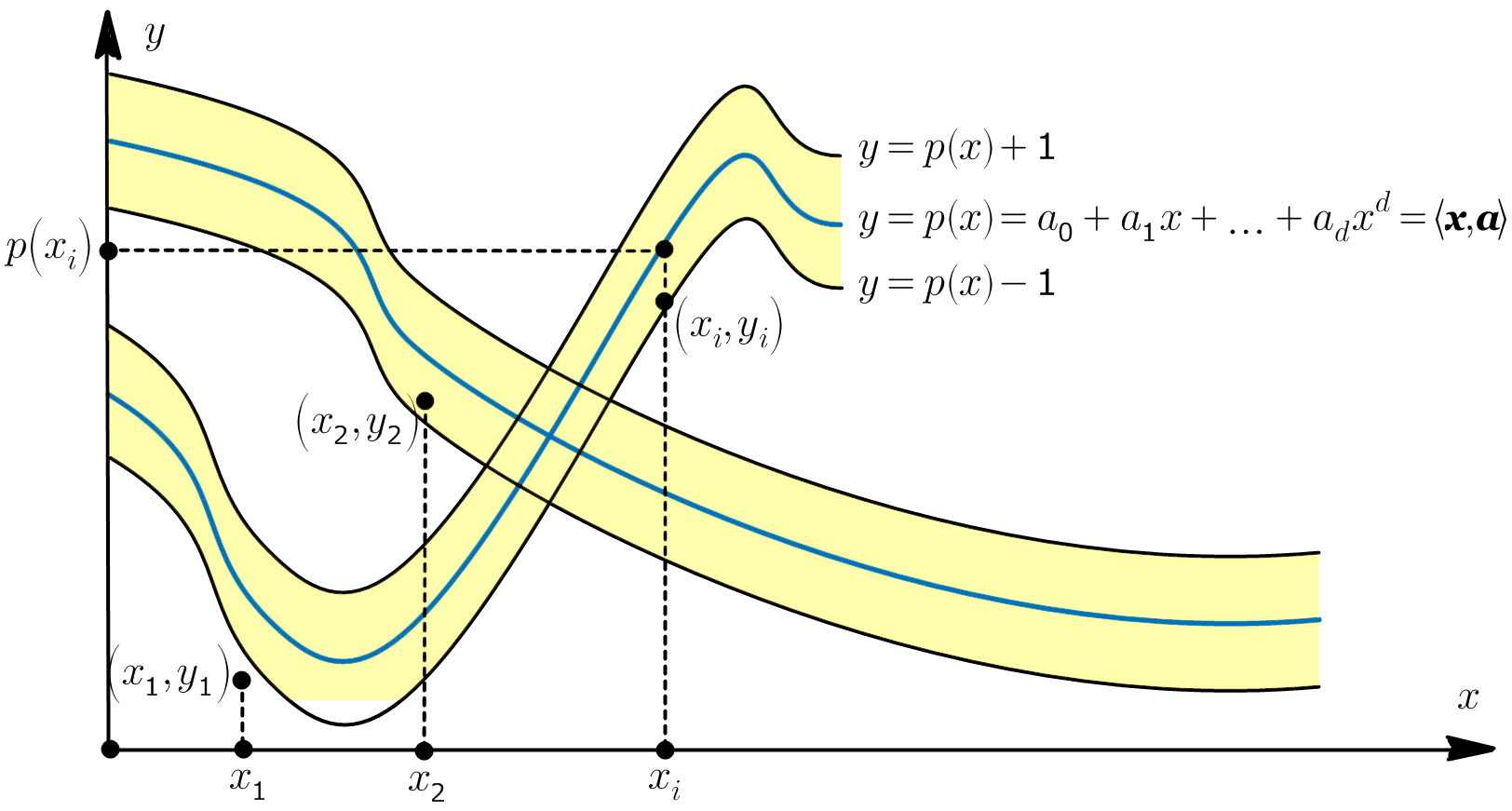}
\smallskip
\caption{Controlling polynomials of degree at most $d$.}
\medskip
\end{figure}

Write $p(x)$ in the form $p(x)=\langle \mathbf{x},\mathbf{a}\rangle,$  where $\mathbf{x} = (1,x,\ldots,x^d)$, $\mathbf{a} = (a_0,a_1,\ldots, a_d)\in{\mathbb{R}}^{d+1}$, and $\langle.\rangle$ stands for the scalar product. Using this notation, we have $\mathbf{x}_i=(1,x_i,\ldots,x_i^d)$ and the inequality $|p(x_i)-y_i|\le 1$ can be rewritten as $$y_i-1\le \langle\mathbf{x}_i, \mathbf{a}\rangle\le y_i+1.$$
For a fixed $i$, the locus of points $\mathbf{a}\in{\mathbb{R}}^{d+1}$ satisfying this double inequality is a slab $S_i\subset{\mathbb{R}}^{d+1}$  of width $w_i = \frac 2{\|\mathbf{x}_i\|} = \frac 2{(\sum_{j=0}^d x_i^{2j})^{1/2}}$, with normal vector $\mathbf{x}_i$. The sequence $x_1, x_2,\ldots$ is ${\mathcal P}_d$-controlling if and only if the sequence of slabs $S_1, S_2,\ldots$ permits a translative covering of $\R^{d+1}$.

\medskip

We distinguish two cases. If $x_i\le3$ for infinitely many (and, hence, for all) integers $i$, then for the widths of the corresponding slabs we have $w_i>\frac 1{3^{d}}$. Thus, these slabs permit a translative covering of $\R^{d+1}$.

\smallskip

Suppose next that $x_i>3$ for all $i\ge m$. Then the sequence $w_{m}\ge w_{m+1}\ge w_{m+2}\ge\ldots$ satisfies the condition of Theorem \ref{th1}. Indeed, we have $w_i \ge \frac 1 {x_i^d}$ for all $i\ge m$, which implies that
$$\limsup_{n\rightarrow\infty}\frac{w_{m}+w_{m+1}+\ldots+w_n}{\log(1/w_n)}\ge
 \limsup_{n\rightarrow\infty}\left(\frac1{x_1^{d}}+\frac1{x_2^{d}}+\ldots+\frac1{x_n^{d}}\right)\Big/{\log x_n}>0$$
Thus, by Theorem~\ref{th1}, the sequence of slabs $S_{m}, S_{m+1}, S_{m+2},\ldots$ permits a translative covering of $\R^{d+1}$, which in turn implies that the sequence $x_1\le x_2\le\ldots$ is ${\mathcal P}_d$-controlling. \qed

\medskip

Note that the above proof yields that for a sequence of positive numbers $x_1\le x_2\le \ldots$ to be ${\mathcal P}_d$-controlling, it is {\em necessary} that $\sum_{i=1}^{\infty}\frac2{x_i^d}\ge\sum_{i=1}^{\infty}w_i=\infty$. Moreover, if the translative plank conjecture stated in the previous section is true, this condition is also {\em sufficient}.
\medskip

In the proof of Corollary~\ref{th3}, we applied Theorem~\ref{th1} to a very special sequence of slabs $S_i$, whose normal vectors lie on a moment curve $\gamma(x)=(1,x,x^2,\ldots,x^d)\subset\mathbb{R}^{d+1}$. Exploring the natural ordering of these vectors along the curve, in Section 5 we will be able to show that the above condition is indeed necessary and sufficient, without proving the translative plank conjecture.

\begin{thm}\label{th31} Let $d$ be a positive integer and $x_1\le x_2\le \ldots$ be a monotone increasing infinite sequence of positive numbers. The sequence $x_1, x_2,\ldots$ is {\em ${\mathcal P}_d$-controlling} if and only if
$\sum_{i=1}^{\infty}\frac1{x_i^d}=\infty$.
\end{thm}

This theorem settles Conjecture 3.2.B in~\cite{MP}.

\section{Proofs of Theorems~\ref{th1} and~\ref{th2}}

We start with the proof of Theorem~\ref{th2}. Then we deduce Theorem~\ref{th1} from Theorem~\ref{th2}.

The proof of Theorem~\ref{th2} is based on some ideas that go back (at least) to {\em Claude Ambrose Rogers}~\cite{Ro,ER}.

\medskip
\noindent{\bf Proof of Theorem~\ref{th2}.} Every slab $S\subset\mathbb{R}^d$ can be expressed in the form
$$S = \{\mathbf{x}\in \mathbb{R}^d: \bar{b}\le \langle {\bf v},\mathbf{x}\rangle\le \bar{b}+w \},$$
where ${\bf v}$ and $w$ are the unit normal vector and the width of $S$, respectively, and $\bar{b}$ is a suitable real number.
\medskip

Fix a sequence of slabs $S_1,\ldots,S_n\subset{\mathbb{R}}^d$ meeting the requirements of the theorem. We may assume that $w_1<1$, otherwise we can cover a ball $B$ of diameter $1$ with a translate of the first slab. Consider the modified sequence of slabs $S'_1,\ldots,S'_n$, where each $S'_i$ is obtained from $S_i$ by reducing its width by a factor of $2$. More precisely, if $$S_i = \{\mathbf{x}\in \mathbb{R}^d: \bar{b}_i\le \langle {\bf v}_i,\mathbf{x}\rangle\le \bar{b}_i+w_i \},$$ then let $$S'_i = \{\mathbf{x}\in \R^d: \bar{b}_i+w_i/4\le \langle {\bf v}_i,\mathbf{x}\rangle\le \bar{b}_i+3w_i/4 \}.$$

\medskip

We describe a greedy algorithm to cover a large part of the unit diameter ball $B$ with suitable translates of $S'_1,\ldots, S'_n$. Set $K_0=B$. For $i=1,\ldots, n$, let $K_{i-1}$ denote the set of points of $B$ not covered by the translates of $S'_1,\ldots,S'_{i-1}$ selected during the first $i-1$ steps of the algorithm. In step $i$, we choose a constant $b_i$ so that the translate $$T(S'_i)=\{\mathbf{x}\in \R^d: b_i+w_i/4\le \langle {\bf v}_i,\mathbf{x}\rangle\le b_i+3w_i/4 \}$$ covers at least a $w_i/3$-fraction of the volume of $K_{i-1}$, i.e., we have
$$Vol(K_{i-1}\backslash T(S'_i))\le (1-w_i/3)Vol(K_{i-1}).$$ Then $K_i=K_{i-1}\backslash T(S'_i)$. It is not difficult to see that such a $b_i$ exists. Indeed, fix a sequence $b_i^1,\ldots, b_i^{\lceil 2/w_i\rceil}$ with the property that the union of translates of $S'_i$,
$$\bigcup_{j=1}^{\lceil 2/w_i\rceil}\{\mathbf{x}\in \R^d: b_i^j+w_i/4\le \langle {\bf v_i},\mathbf{x}\rangle\le b_i^j+3w_i/4 \},$$
completely covers $B$. Since $K_{i-1}\subseteq B$, it follows by the pigeonhole principle that at least one of these translates will cover at least a $\frac 1{\lceil 2/w_i\rceil}$-fraction of the volume of $K_{i-1}$. Notice that $\frac 1{\lceil 2/w_i\rceil}>\frac 1{1+2/w_i} = \frac {w_i}{2+w_i}>\frac{w_i}3.$
\medskip

\begin{figure}[h]
\medskip
\centering
\smallskip
\includegraphics[height=3in]{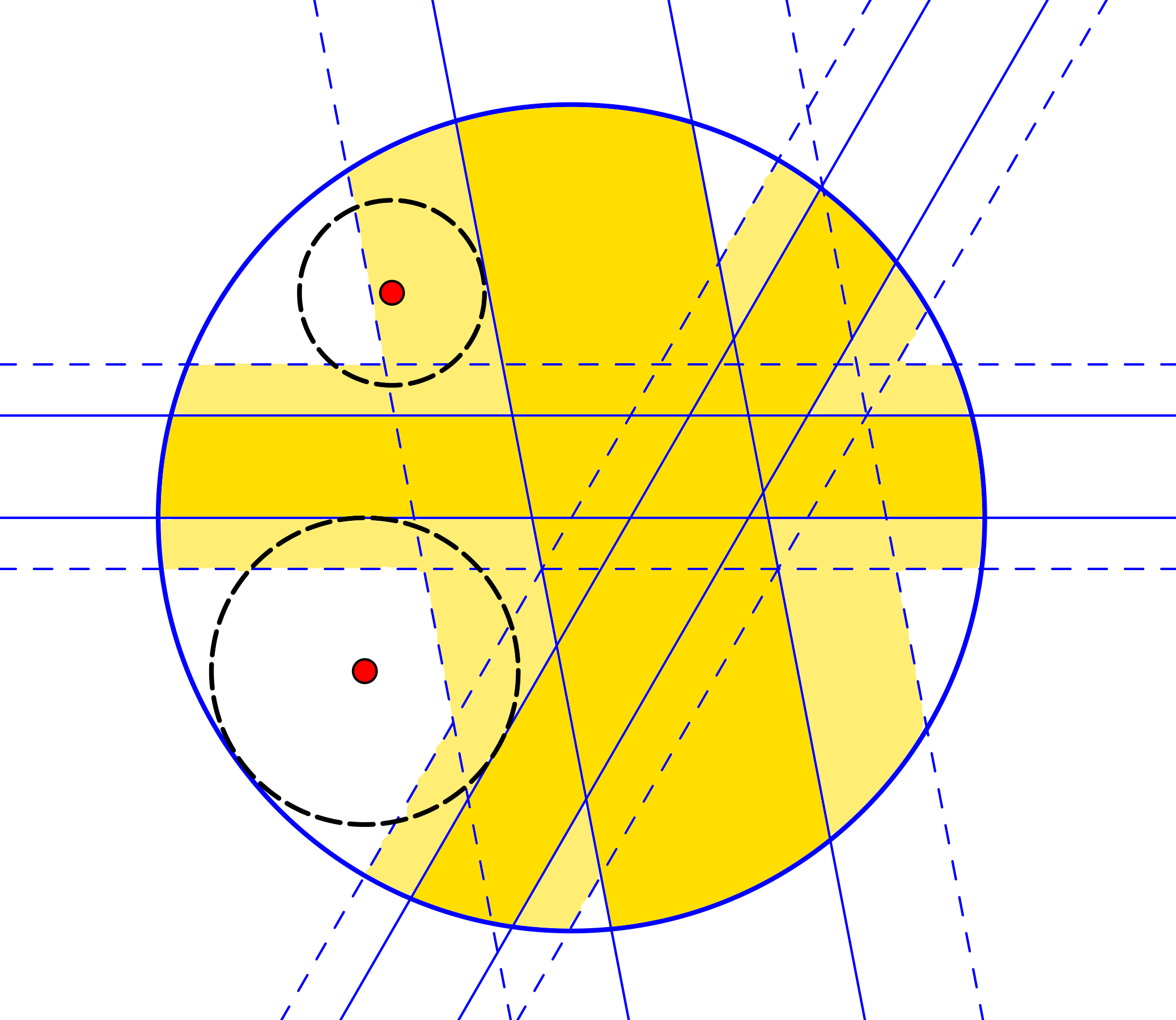}
\smallskip
\caption{After $n$ steps, the uncovered part of $B$ contains no ball of radius $w_n/4$.}
\medskip
\end{figure}

After $n$ steps, the volume of the set of uncovered points of $B$ satisfies
$$Vol(K_n)\le Vol(B)\prod_{i=1}^n(1-w_i/3)< Vol(B)\exp\left\{-\frac13\sum_{i=1}^n w_i\right\}.$$
Using the assumption on $\sum_{i=1}^n w_i$, we obtain
$$Vol(K_n)< Vol(B)\exp\{-d\log(2/w_n)\}=Vol(B)\left(\frac {w_n}2\right)^d.$$

\noindent Therefore, the set $K_n\subseteq B$ does not contain a ball of radius $w_n/4$. See Fig. 3. This implies that every point of $B$ is at distance at most $w_n/4$ from the surface of $B$ or from one of the selected translates $T(S'_i)\; (i=1,\ldots,n)$. In other words, expanding each $T(S_i')$ by a factor of $2$ around its hyperplane of symmetry, we obtain a translate of $S_i$, and the union of these translates,
$$\bigcup_{i=1}^{n}\{\mathbf{x}\in \R^d: b_i\le \langle {\bf v}_i,\mathbf{x}\rangle\le b_i+w_i \},$$
covers the ball of radius $1/2-w_n/4$, concentric with $B$. \qed

\bigskip
\noindent{\bf Proof of Theorem~\ref{th1}.}
It follows from the condition of the theorem that the sequence of slabs $S_i$ can be split into infinitely many finite subsequences ${\mathcal{S}}_j=(S_{i_j+1}, S_{i_j+2},\ldots, S_{i_{j+1}})$ for $j=1, 2,\ldots,$ where $0=i_1 <i_2<i_3<\ldots,$  and
\begin{equation}\label{eq1} w_{i_j+1}+ w_{i_j+2}+\ldots + w_{i_{j+1}}\ge  c\log(1/w_{i_{j+1}}),\end{equation}
for a suitable positive constant $c\le 1$.

We can assume without loss of generality that $w_i\le \frac{c}{3d}$ holds for every $i$. Otherwise, if there are finitely many exceptional indices $i$ with $w_i>\frac{c}{3d}$, we simply discard the corresponding slabs. The remaining sequence will meet the requirements of Theorem~\ref{th1}. If the number of exceptional indices is infinite, then we use the construction of the covering for sequences of slabs with bounded widths, described in the beginning of Section 2.
\medskip

For the same reason, it suf\-fices to show that each sub\-se\-quence
${\mathcal{S}}_j=(S_{i_j+1}, S_{i_j+2},\ldots, S_{i_{j+1}})$ permits a translative covering of a ball of diameter $\frac{c}{6d}$. It follows from Theorem~\ref{th2} that if
\begin{equation}\label{eq2} \frac{3d}c(w_{i_j+1}+ w_{i_j+2}+\ldots + w_{i_{j+1}})\ge
3d \log(2/(\frac{3d}{c}w_{i_{j+1}})),\end{equation}
then the slabs $\frac{3d}cS$ obtained from the elements $S\in {\mathcal{S}}_j$ by widening them by a factor of $\frac{3d}c$ permit a translative covering of a ball of diameter $1-\frac{3d}{2c}w_{i_{j+1}}$. Therefore, by scaling, the strips $S_i$ permit a translative covering of a ball of the diameter $\frac{c}{3d}(1-\frac{3d}{2c}w_{i_{j+1}})\ge\frac{c}{6d}.$
In the last step, we used the assumption $w_{i_{j+1}}\le \frac{c}{3d}$.

\medskip

It remains to verify is that (\ref{eq1}) implies (\ref{eq2}), but this reduces to the inequality $2c\le3d$. \qed

\section{Proof of Theorem~\ref{th31}}

It follows from the remark right after the proof of Corollary~\ref{th3} that we only have to establish the ``if'' part of the theorem. That is, using the same notation as in the proof of Corollary~\ref{th3}, it is sufficient to show that the sequence of slabs $S_i$ defined there permits a translative covering, provided that the sum $\frac1{x_{1}^{d}}+\frac1{x_{2}^{d}}+\ldots$ is divergent. As we have seen, we can also assume without loss of generality that $x_i\ge 3$ for every $i$.

\medskip

In the proof of Corollary~\ref{th3}, we blindly applied Theorem~\ref{th1} to the slabs $S_i$, without exploring their special properties. To refine our argument, we are going to exploit the specifics of the slabs. Recall that $S_i$ has width $w_i >\frac1{x_i^{d}}$ for every $i$, and its normal vector $\mathbf{x}_i = (1,x_i,\ldots,x_i^d)$ lies on the moment curve $(1,x,x^2,\ldots,x^d)$. First, we need an auxiliary lemma.

\begin{lem}\label{lem2} Let $d$ be a positive integer, let $3\le x_1\le x_2\le\ldots$ be a finite or infinite sequence of reals, and let $\mathbf{x}_i=(1,x_i,x_i^2,\ldots,x_i^d)$ for every $i$.
Then there exist $d+1$ linearly independent vectors $\mathbf{u}_1,\ldots,\mathbf{u}_{d+1}\in \mathbb{R}^{d+1}$ such that for every $i\; (i=1,2,\ldots)$ and $j\; (j=1,2,\ldots,d+1)$, we have

$$(i)\;\;\;\;\;\; \frac{\langle\mathbf{x}_{i+1},\mathbf{u}_1\rangle}{\langle\mathbf{x}_{i},\mathbf{u}_1\rangle}\le
\frac{\langle\mathbf{x}_{i+1},\mathbf{u}_j\rangle}{\langle\mathbf{x}_{i},\mathbf{u}_j\rangle},$$

$$(ii) \;\;\;\;\;\;\
\langle\mathbf{x}_i,\mathbf{u}_j\rangle\ge \frac13\|\mathbf{x}_i\|\|\mathbf{u}_j\|.$$
\end{lem}

\bigskip
\noindent{\bf Proof.} Take the standard basis $\mathbf{e}_1,\ldots,\mathbf{e}_{d+1}$ in $\R^{d+1}$, i.e., let $e_i$ denote the all-zero vector with a single $1$ at the $i$-th position. Set $\mathbf{u}_j := \mathbf{e}_{d+1-j}+\mathbf{e}_{d+1}$ for $j=1,\ldots, d$ and $\mathbf{u}_{d+1} :=\mathbf{e}_{d+1}$.

Condition (i) trivially holds for $j=1$ and it is very easy to check it for $j=d+1$. For $j=2,\ldots,d,$ it reduces to
$$\frac{x_{i+1}^{d-1}+x_{i+1}^{d}}{x_{i}^{d-1}+x_{i}^{d}}\le \frac{x_{i+1}^{d-j}+x_{i+1}^{d}}{x_{i}^{d-j}+x_{i}^{d}},$$
which is equivalent to
$$(x_{i+1}^{d-1}+x_{i+1}^{d})(x_{i}^{d-j}+x_{i}^{d})\le (x_{i+1}^{d-j}+x_{i+1}^{d})(x_{i}^{d-1}+x_{i}^{d}).$$
The last inequality can be rewritten as
$$x_{i+1}^{d-j}x_{i}^{d-j}(x_{i+1}-x_i)(\sum_{k=0}^{j-1}x_{i+1}^kx_i^{j-1-k}+\sum_{k=0}^{j-2}x_{i+1}^kx_i^{j-2-k} - x_{i+1}^{j-1}x_{i}^{j-1})\le 0,$$
or, dividing both sides by $x_{i+1}^{d-j}x_{i}^{d-j}(x_{i+1}-x_i)$, as
$$\sum_{k=0}^{j-1}x_{i+1}^kx_i^{j-1-k}+\sum_{k=0}^{j-2}x_{i+1}^kx_i^{j-2-k}-x_{i+1}^{j-1}x_{i}^{j-1}\le 0.$$

Using the fact $x_{i+1}\ge x_{i}$, and bounding from above each sum by its largest term multiplied by the number of terms, we obtain that the left-hand side of the last inequality is at most
$$jx_{i+1}^{j-1}+(j-1)x_{i+1}^{j-2}- x_{i+1}^{j-1}x_{i}^{j-1}<
x_{i+1}^{j-1}(2j-1-x_i^{j-1}).$$
As $x_i\ge 3$, the right-hand side of this inequality is always nonpositive and (ii) holds.

\smallskip

It remains to verify condition (ii). Taking into account that $x_i\ge 3$, we have
$$\langle \mathbf{x}_i,\mathbf{u}_{d+1}\rangle = x_i^{d}\ge \frac12\|\mathbf{x}_i\| = \frac12\|\mathbf{x}_i\|\|\mathbf{u}_{d+1}\|.$$
On the other hand, for $j = 1,\ldots,d$, we obtain
$$\langle\mathbf{x}_i, \mathbf{u}_{j}\rangle = x_i^{d-j}+x_i^{d}\ge \frac12\|\mathbf{x}_i\| \ge \frac13\|\mathbf{x}_i\|\|\mathbf{u}_j\|.$$
This completes the proof of Lemma \ref{lem2}. \qed

\bigskip

In order to establish Theorem~\ref{th31}, it is enough to prove that there is a constant $c=c(d+1)$ such that any system of slabs $S_i\; (i=1,\ldots,n)$ in $\mathbb{R}^{d+1}$ whose normal vectors are $(1,x_i,\ldots,x_i^d)$ for some $3\le x_1\le x_2\le\ldots\le x_n$ and whose total width is at least $c$, permits a translative covering of a ball of unit diameter. This is an immediate corollary of Lemma~\ref{lem2} and the following assertion.

\begin{lem}\label{lem1}
For every positive integer $d$, for any system of $d+1$ linearly independent vectors $\mathbf{u}_1,\ldots,\mathbf{u}_{d+1}$ in $\mathbb{R}^{d+1}$, and for any $\gamma>0$, there is a constant $c$ with the following property.

Given any system of slabs $S_i\; (i=1,\ldots,n)$ in $\mathbb{R}^{d+1}$, whose normal vectors $\mathbf{x}_i$ satisfy the conditions

$$(i) \;\;\;\;\;\;
\frac{\langle\mathbf{x}_{i+1},\mathbf{u}_1\rangle}{\langle\mathbf{x}_{i},\mathbf{u}_1\rangle}\le
\frac{\langle\mathbf{x}_{i+1},\mathbf{u}_j\rangle}{\langle\mathbf{x}_{i},\mathbf{u}_j\rangle},$$

$$(ii) \;\;\;\;\;\;\;\;
\langle\mathbf{x}_i,\mathbf{u}_j\rangle\ge \gamma\|\mathbf{x}_i\|\|\mathbf{u}_j\|$$
\smallskip

\noindent for every $i$ and $j$, and whose total width $\sum_{i=1}^nw_i$ is at least $c$, the slabs $S_i$ permit a translative covering of a $(d+1)$-dimensional ball of unit diameter.
\end{lem}

The proof is based on a greedy algorithm for covering a large simplex in $\mathbb R^{d+1}$ by the translates of $S_i$. One of the vertices of the covered simplex is $\mathbf 0$, and the others lie on the rays emanating from $\mathbf 0$ in the directions $\mathbf u_1,...,\mathbf u_{d+1}$. At each step, we place a new slab in such a way that the newly covered part of the ray in direction $\mathbf u_1$ adjoins the previously covered part. Condition (i) guarantees that we do not leave any ``hole" between the translates uncovered. Condition (ii) ensures that the simplex completely covered at the end of the procedure is "non-degenerate": the ratio of any two of its sides is bounded from above by an absolute constant. Finally, we will show that the side of this simplex along $\mathbf u_1$, and thus every other side of it, is sufficiently large, which implies that the simplex contains a unit ball.
\bigskip

\noindent{\bf Proof.}
Instead of covering a ball of unit diameter, it will be more convenient to cover the simplex $\Delta$ with one vertex in the origin $\mathbf{0}$ and the others at the points (vectors) $\mathbf{u}_j\; (j=1,\ldots,d+1)$. By properly scaling these vectors, if necessary, we can assume that $\Delta$ contains a ball of unit diameter.

\begin{figure}[h]
\centering
\smallskip
\includegraphics[height=3in]{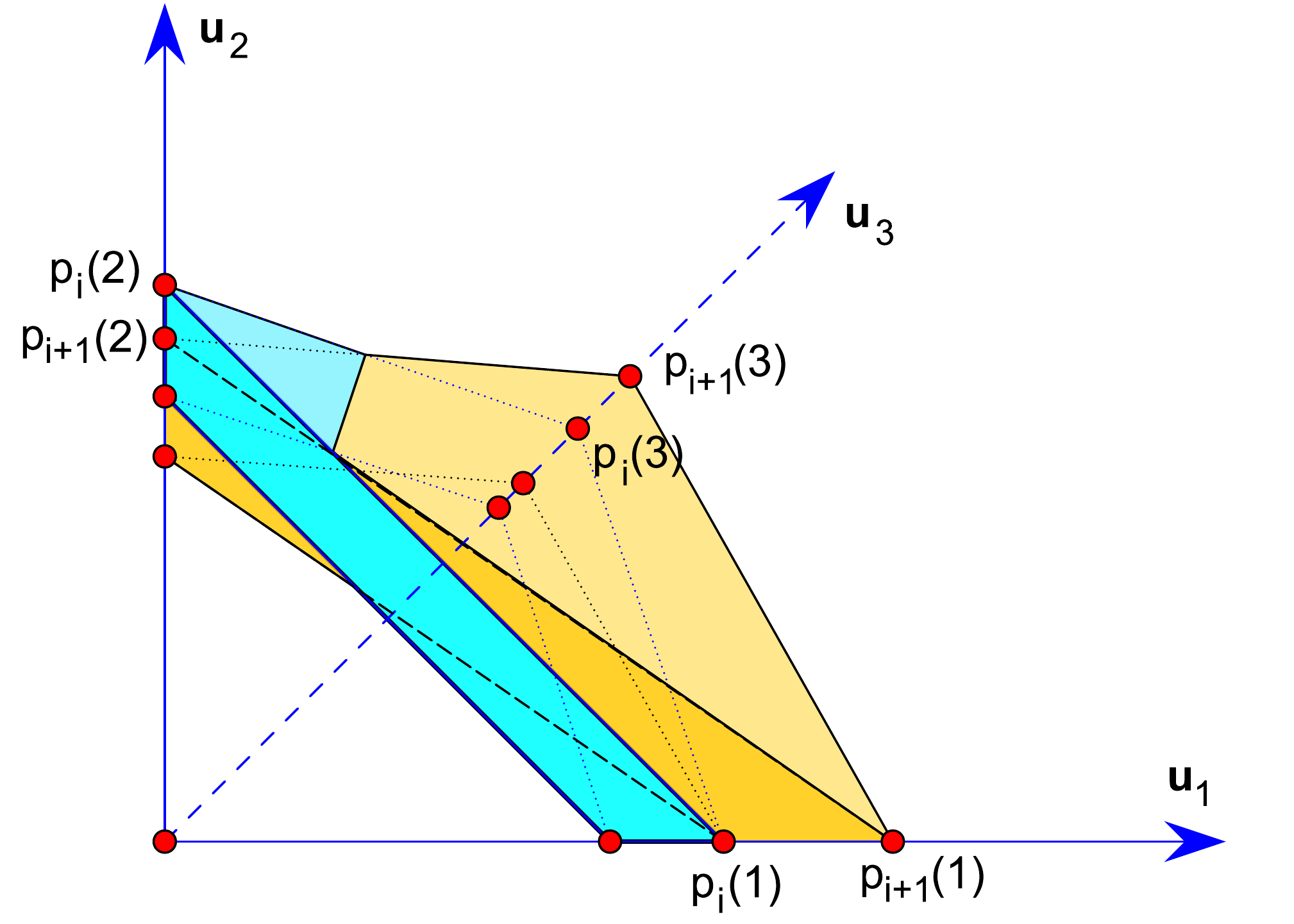}
\smallskip
\caption{We place the slabs one by one.}
\medskip
\end{figure}

\medskip

We place the slabs one by one. See Fig. 4. We place $S'_1$, a translate of $S_1$, so that one of its boundary hyperplanes passes through $\mathbf{0}$ and the other one cuts a simplex $\Delta_1$ out of the cone $\Gamma$ of all linear combinations of the vectors $\mathbf{u}_1,\ldots,\mathbf{u}_{d+1}$ with positive coefficients. According to our assumptions, we have $\langle\mathbf{x}_1,\mathbf{u}_j\rangle>0$ for every $j$. Therefore, $S'_1$ does not separate $\Gamma$ into two cones: $S'_1\cap\Gamma$ is indeed a simplex $\Delta_1$.

\medskip

We place a translate $S_2'$ of $S_2$ in such a way that it adjoins $\Delta_1$ in its vertex along the ray emanating from $\mathbf 0$ in direction $\mathbf u_1$. As we will see later, due to condition (i), $S_2'$ and $\Delta_1$ overlap (in at least one point) along the rays in direction $\mathbf u_i$ for $i\ge 2$.  Suppose that we have already placed $S'_1,\ldots,S'_i$, the translates of $S_1,\ldots,S_i,$ so that their union covers a simplex $\Delta_i$ with one vertex at the origin, and the others along the $d+1$ half-lines that span the cone $\Gamma$. We also assume that the facet of $\Delta_i$ opposite to the origin is a boundary hyperplane of $S'_i$. Let $\mathbf{p}_i(j)$ denote the vertex of $\Delta_i$ that belongs to the open half-line parallel to $\mathbf{u}_j$ emanating from $0\; (j=1,\ldots,d+1)$.
\medskip

Next, we place a translate $S'_{i+1}$ of $S_{i+1}$ so that one of its boundary hyperplanes, denoted by $\pi$, passes through $\mathbf{p}_i(1)$, and the other one, $\pi'$, cuts the half-line parallel to $\mathbf{u}_1$ at a point $\mathbf{p}_{i+1}(1)$ with $\|\mathbf{p}_{i+1}(1)\|>\|\mathbf{p}_i(1)\|$. That is, $\mathbf{p}_{i+1}(1)$ is further away from the origin than $\mathbf{p}_i(1)$ is. Let $\mathbf{p}_{i+1}(2),\ldots,\mathbf{p}_{i+1}(d+1)$ denote the intersection points of $\pi'$ with the half-lines parallel to $\mathbf{u}_2,\ldots,\mathbf{u}_{d+1}$, respectively, and let $\Delta_{i+1}$ be the simplex induced by the vertices $\mathbf{0},\mathbf{p}_{i+1}(1),\ldots,\mathbf{p}_{i+1}(d+1)$.

\medskip

We have to verify that $\Delta_{i+1}$ is entirely covered by the slabs $S'_1,\ldots, S'_{i+1}$. By the induction hypothesis, $\Delta_i$ was covered by the slabs $S'_1,\ldots,S'_i.$ Thus, it is sufficient to check that the hyperplane $\pi$ intersects every edge $\mathbf{0}\mathbf{p}_i(j)$ of $\Delta_{i}$, for $j=1,\ldots,d+1$. Let $\alpha_j\mathbf{u}_j$ be the intersection point of $\pi$ with the half-line parallel to $\mathbf{u}_j$, and let $\mathbf{p}_i(j)=\beta_j\mathbf{u}_j$. We have to prove that $\alpha_j\le\beta_j$.

\medskip

By definition, we have $\langle \mathbf{x}_{i+1}, \mathbf{p}_i(1)-\alpha_j\mathbf{u}_j\rangle=0$ and $\langle\mathbf{x}_i,
\mathbf{p}_i(1)-\beta_j\mathbf{u}_j\rangle=0$. From here, we get
$$\frac{\alpha_j}{\beta_j}=
\frac{\langle\mathbf{x}_{i+1},\mathbf{p}_i(1)\rangle}{\langle\mathbf{x}_{i+1},\mathbf{u}_j\rangle}
\Big/\frac{\langle\mathbf{x}_{i},\mathbf{p}_i(1)\rangle}{\langle\mathbf{x}_{i},\mathbf{u}_j\rangle}
=\frac{\langle\mathbf{x}_{i+1},\mathbf{p}_i(1)\rangle}{\langle\mathbf{x}_{i},\mathbf{p}_i(1)\rangle}
\Big/\frac{\langle\mathbf{x}_{i+1},\mathbf{u}_j\rangle}{\langle\mathbf{x}_{i},\mathbf{u}_j\rangle}
=\frac{\langle\mathbf{x}_{i+1},\mathbf{u}_1\rangle}{\langle\mathbf{x}_{i},\mathbf{u}_1\rangle}
\Big/\frac{\langle\mathbf{x}_{i+1},\mathbf{u}_j\rangle}{\langle\mathbf{x}_{i},\mathbf{u}_j\rangle}.$$
In view of assumption (i) of the lemma, the right-hand side of the above chain of equations is at most $1$, as required.

\medskip

Observe that during the whole procedure the uncovered part of the cone $\Gamma$ always remains convex and, hence, connected. In the $n$th step, $\cup_{i=1}^nS'_i\supset\Delta_n$. By the construction, $\mathbf{p}_i(1)$ lies at least $w_i$ farther away from the origin along the half-line parallel to $\mathbf{u}_1$ than $\mathbf{p}_{i-1}(1)$ does. Thus, we have
$$\|\mathbf{p_n}(1)\|\ge \sum_{i=1}^nw_i\ge c.$$

Using the fact that $\langle\mathbf{x}_n,\mathbf{p}_n(j)-\mathbf{p}_n(1)\rangle=0$ for every $j\ge 2$, and taking into account assumption (ii), we obtain
$$\|\mathbf{p}_n(j)\|\ge\frac{\langle\mathbf{x}_n,\mathbf{p}_n(j)\rangle}{\|\mathbf{x}_n\|}
=\frac{\langle\mathbf{x}_n,\mathbf{p}_n(1)\rangle}{\|\mathbf{x}_n\|}
\ge\gamma\|\mathbf{p_n}(1)\|\ge\gamma c.$$
Thus, if $c$ is sufficiently large, we have $\|\mathbf{p}_n(j)\|\ge\|\mathbf{u}_j\|$. This means that $\Delta_n$ contains the simplex $\Delta$ defined in the first paragraph of this proof. Hence, it also contains a ball of unit diameter, as required. \qed

\section{Concluding remarks}

\noindent{\bf 1.} As was mentioned in Section 2, the translative packing conjecture of Makai and Pach is known to be true in the plane. Moreover, in~\cite{MP} a stronger statement was proved: there exists a constant $c$ such that every collection of strips with total width at least $c$ permits a translative covering of a disk of diameter $1$. In view of this, one can make the following even bolder conjecture.
\medskip

\noindent{\bf Strong translative plank conjecture.} {\em For any positive integer $d$, there exists a constant $c=c(d)$ such that every sequence of slabs in $\R^d$ with total width at least $c$ permits a translative covering of a unit diameter $d$-dimensional ball.}
\medskip

Suppose that the {\em translative plank conjecture} (see Section 2) is true for a positive integer $d$. Answering a question in \cite{MP}, {\em Imre Z. Ruzsa}~\cite{Ru} proved that then, for the same value of $d$, the {\em strong translative plank conjecture} also holds. Thus, the two conjectures are equivalent.
\medskip

\noindent{\bf 2.} We say that a sequence of positive numbers $x_1\le x_2\le\ldots$ is {\em strongly ${\mathcal P}_d$-controlling} if there exist reals $y_1, y_2,\ldots$ with the property that, for every $\varepsilon>0$ and for every polynomial $p$ of degree at most $d$, one can find an $i$ with $$|f(x_i)-y_i|\le \varepsilon.$$

It is easy to see that the condition in Theorem~\ref{th31} is sufficient to guarantee that the sequence $x_1,x_2,\ldots$ strongly controls ${\mathcal P}_d$ for every $d$.

Controlling sequences can be analogously defined for any other class of functions $f:\mathbb{R}^k\rightarrow\mathbb{R}^l$. Several other problems of this kind are discussed in~\cite{MP} and~\cite{KP}.

\end{document}